\documentclass[aps,prl,preprint,groupedaddress,amsfonts,floatfix]{revtex4-1}

\usepackage{longtable}
\usepackage{color}
\usepackage{graphicx}
\usepackage{epstopdf}
\usepackage{amsmath}

\makeatletter
\renewcommand{\@biblabel}[1]{\quad#1.}

\makeatother

\begin{document}


\title{Continuous-time statistics and generalized relaxation equations}


\author{Enrico Scalas}
\email[]{e.scalas@sussex.ac.uk}
\affiliation{Department of Mathematics, School of Mathematical and Physical Sciences, University of Sussex, Brighton, UK}


\date{\today}

\begin{abstract}
Using two simple examples, the continuous-time random walk as well as a two state Markov chain, the relation between generalized anomalous relaxation equations and semi-Markov processes is illustrated. This relation is then used to discuss continuous-time random statistics in a general setting, for statistics of convolution-type. Two examples are presented in some detail: the sum statistic and the maximum statistic.
\end{abstract}

\maketitle
\section{Introduction}
Continuous-time random walks (CTRWs) are a straightforward generalization of compound Poisson processes. Their simplest version, the so-called {\em uncoupled} case, can be defined as follows. Let $\{Y_i\}_{i=1}^\infty$ be a sequence of independent and identically distributed random variables in $\mathbb{R}^d$ (here, for the sake of simplicity, we consider $d=1$) with cumulative distribution function $F_{Y_1}(u)=\mathbb{P}(Y_1 \leq u)$. The corresponding {\em random walk} is the homogeneous Markov chain defined by
\begin{equation}
\label{randomwalk}
X_n = \sum_{i=1}^n Y_i.
\end{equation} 
Now, suppose we are given a sequence of positive independent and identically distributed random variables $\{J_i \}_{i=1}^\infty$ with the meaning of inter-event durations and with cumulative distribution function $F_{J_1} (w) = \mathbb{P}(J_1 \leq w)$. Further assume that the sequences $\{Y_i\}_{i=1}^\infty$ and $\{J_i \}_{i=1}^\infty$ are independent. First define the {\em epochs} at which events occur as
\begin{equation}
\label{epochs}
T_n = \sum_{i=1}^n J_i,
\end{equation}
then introduce the number of events from $T_0 = 0$ seen as an event (technically, as a {\em renewal point})
\begin{equation}
\label{countingprocess}
N(t) = \max \{n: T_n \leq t\}.
\end{equation} 
Change time from $n$ to $N(t)$ in eq. \eqref{randomwalk} to get the uncoupled continuous-time random walk
\begin{equation}
\label{CTRWuncoupled}
X(t) = X_{N(t)} = \sum_{i=1}^{N(t)} Y_i.
\end{equation}

If the $J_i$s are exponentially distributed, then $N(t)$ is a Poisson process and eq. \eqref{CTRWuncoupled} defines compound Poisson processes \cite{feller66,definetti75}. These are L\'evy processes \cite{applebaum09} with L\'evy triplet given by $(0,0,\lambda \sigma)$ where drift and diffusion are $0$ and $\sigma$ is a measure on $\mathbb{R}$ with $\sigma\{0\}=0$. Just as a reminder, a L\'evy process is a Markov process with independent and stationary increments. The realizations (a.k.a sample paths) of a L\'evy process are right-continuous with left limits (or {\em c\`adl\`ag} from the French {\em continu \`a droit, limite \`a gauche}). Compound Poisson processes play an important role in the theory of L\'evy processes as they can approximate any other L\'evy process. To be more precise, one can consider any L\'evy process as an independent sum of a Brownian motion with drift and a countable number of independent compound Poisson processes with different jump rates $\lambda$ and jump distributions $\sigma$ \cite{kyprianou}.

If the $J_i$s are not exponentially distributed, then $N(t)$ is a counting renewal process and eq. \eqref{CTRWuncoupled} defines a renewal-reward process that is non-Markovian and non-L\'evy, but semi-Markov \cite{germano09}. The realizations are assumed to be {\em c\`adl\`ag} as well; this is useful for functional limit theorems \cite{meerschaert04}. In \cite{germano09}, a derivation of the famous Montroll-Weiss equation \cite{montroll65} as a necessary condition for semi-Markov processes is presented. In this case, the continuous-time random walk becomes a process with infinite memory. This is due to the infinite memory of the counting process $N(t)$ that, in its turn, is due to the infinite memory of the residual time to the next renewal from any ``observation'' time $t$ \cite{politi11}.

Assume absolute continuity of the distributions of $X_i$s and $J_i$s and define their respective probability density functions as $f_{Y_1} (u) = d F_{Y_1} (u)/du$ for the jumps and $f_{J_1} (w) = d F_{J_1} (w)/dw$ for the inter-event times. Then, straightforward calculations (see the Appendix) lead from the equation of Montroll and Weiss to the following evolution equation \cite{mainardi00}:
\begin{equation}
\label{CTRWevolution}
\int_0^t \Phi(t-t') \frac{\partial}{\partial t'} p(x,t') \, dt'= - p(x,t) + \int_{-\infty}^{+ \infty} f_{Y_1} (x - x') p(x',t) \, dx',
\end{equation}
where $p(x,t) = dF_{X(t)} (x)/dx$ is the probability density function of finding the continuous time random in $x$ at time $t$ given that $X(0) = 0$ and $\Phi(t)$ has the following Laplace transform
\begin{equation}
\label{memorykernel}
\mathcal{L}(\Phi(t))(s) = \frac{1 - \mathcal{L}(f_{J_1} (t))(s)}{s \mathcal{L}(f_{J_1} (t))(s)}.
\end{equation}
It is interesting to remark that eq. \eqref{CTRWevolution} highlights the infinite memory of the process as $\Phi(t)$ plays the role of {\em memory kernel}. The reader might be interested in comparing this approach to semi-Markov processes with the classical approach in \cite{haenggi77} for non-Markovian processes. In the exponential/Poisson case (set $\lambda=1$, for the sake of simplicity), one has $f_{J_1} (t) = \exp(-t)$ and one gets $\mathcal{L}(\exp(-t))(s) = 1/(1+s)$ so that $\mathcal{L}(\Phi(t))(s) =1$ and $\Phi(t) = \delta(t)$; then, eq. \eqref{CTRWevolution} reduces to the so-called {\em Kolmogorov-Feller} equation \cite{saichev97}
\begin{equation}
\label{CTRWevolutionMarkov}
\frac{\partial}{\partial t} p(x,t) \, dt= - p(x,t) + \int_{-\infty}^{+ \infty} f_{Y_1} (x - x') p(x',t) \, dx'.
\end{equation}
\textcolor{black}{Equation \eqref{CTRWevolution} naturally leads to anomalous diffusion when inter-event times have a power-law distribution with infinite first moment (see \cite{hilfer95} and the references quoted at the end of the next section).}

\section{Anomalous relaxation}

Eq. \eqref{CTRWevolution} is an instance of anomalous relaxation equation. These equations describe governing equations for time-changed Markov processes. They have been recently studied in full generality by Meerschaert and Toaldo \cite{meerschaert15} and independently obtained for a specific semi-Markov random graph dynamics by Georgiou et al. \cite{georgiou15}. In order to illustrate the relationship with relaxation processes more convincingly, we can use the simple example of a homogeneous Markov chain $Y_n$ for $n\leq 0$ with two states $A$ and $B$ and transition probabilities $q_{i,j} = P(X_1=j|X_0=i)$ given by
$q_{A,A} =0$, $q_{A,B} = 1$, $q_{B,A} = 0$ and $q_{B,B}=1$ \cite{jurlewicz00}. This means that if the chain is prepared in state $A$, it will jump to state B at the first step and it will stay there forever. Now define, as above, the new process $Y(t) = Y_{N(t)}$. Then, we have (see \cite{georgiou15} for an explicit derivation and remember that $T_0=0$ is a renewal point):
\begin{equation}
\label{Relax1}
p_{i,j} (t) = \mathbb{P}(Y(t)=j|X(0)=i)=  \bar{F}_{J_1} (t)\delta_{i,j} + \sum_{n=1}^\infty q_{i,j}^{(n)} \mathbb{P} (N(t) = n),
\end{equation}
where $\bar{F}_{J_1} (t) = 1 - F_{J_1} (t)$ is the complementary cumulative distribution function of the inter-event durations. Then we immediately have that $p_{A,A} (t) = \bar{F}_{J_1} (t)$. If $N(t)$ is Poisson with parameter $\lambda =1$, it turns out that $p_{A,A} (t) = \exp(-t)$. In other words the probability of finding the chain in state $A$ decays exponentially as $t$ grows. In this case, $p_{A,A} (t)$ is the solution of the following relaxation Cauchy problem
\begin{equation}
\label{Relax2}
\frac{d p_{A,A} (t)}{dt} = - p_{A,A} (t), \;\;\; p_{A,A} (0) = 1.
\end{equation}
For a general renewal counting process $N(t)$, one gets the following generalized relaxation Cauchy problem instead:
\begin{equation}
\label{Relax3}
\int_0^t \Phi(t-t') \frac{d p_{A,A} (t')}{dt'} dt' = - p_{A,A} (t), \;\;\; p_{A,A} (0) = 1.
\end{equation}

\textcolor{black}{The anomalous relaxation theory outlined above was studied by Scher and Montroll \cite{scher75} in the context of transit-time dispersion in amorphous solid. They 
explicitly assumed a power-law behaviour for the distribution of the inter-event durations. This theory was further developed by Klafter and Silbey \cite{klafter80} who studied transport of particles on a lattice using the projector operator technique. They showed that the exact equation governing the transport averaged over all configurations can be written either as a generalized master equation or as the continuous-time random-walk equations.}

\textcolor{black}{In the next section, the relation will be presented between the anomalous relaxation discussed in \cite{scher75} and fractional operators. This was already discussed in papers by Gl\"ockle and Nonnenmacher and Metzler {\em et al.} \cite{gloeckle91,gloeckle93,metzler95}. Mainardi, Gorenflo and co-workers published two review papers on anomalous relaxation and fractional calculus \cite{mainardi07,oliveira14}. Two general reviews on fractional diffusion, Fokker-Planck equations, and relaxation equations were written by Metzler and Klafter \cite{metzler00,metzler04}. More recently, important properties of CTRWs such as ageing or weak ergodicity breaking have been reviewed as well \cite{sokolov12,metzler14}.}

\textcolor{black}{The useful character of the simple idea of CTRWs fully emerges from the body of work outlined above. Therefore, it is not surprising to see that the time-change from a deterministic $n$ to the random process $N(t)$ can lead to further developments.}

\section{Theory}

After establishing the connection between the random time change $N(t)$ and relaxation equations of the type \eqref{CTRWevolution} and \eqref{Relax3}, one can proceed to study continuous-time statistics in a rather general way.
Let $\{X_i\}_{i=1}^n$ be a sequence of $n$ independent and identically distributed positive random variables with cumulative distribution function $F_{X_1}(u)=\mathbb{P}(X_1 \leq u)$. A {\em statistic}
is a function from $\mathbb{R}^n$ to $\mathbb{R}$ that summarizes some characteristic behavior of the random variables:
\begin{equation}
\label{statistic}
S_n = G_n(X_1, \ldots X_n).
\end{equation}
The statistic $S_n$ is a random variable and, usually, something is known on its distribution. Asymptotic analytical results may be 
available in the limit of large $n$ and 
Monte Carlo simulations can be used for small values of $n$. Let $F_{S_n} (u) = \mathbb{P}(S_n \leq u)$ denote the cumulative 
distribution function of $G_n$. As in the two previous examples, in order to introduce continuous-time statistics, we
use another set of positive independent and identically distributed random variables (independent from the $X_i$s)
$\{J_i\}_{i=1}^\infty$ with the meaning of sojourn times. Let $F_{J_1} (t) = \mathbb{P}(J \leq t)$ denote the cumulative distribution function of the $J_i$s and $f_{J_1} (t) = d F_{J_1} (t)/dt$ denote their probability density function. We again introduce the epochs at which
events occur
\begin{equation}
\label{epochs}
T_n = \sum_{i=1}^n J_i,
\end{equation}
and the counting process $N(t)$ giving the number of events that occur up to time $t$
\begin{equation}
\label{counting}
N(t) = \max \{n: \; T_n \leq t\}.
\end{equation}
The continuous-time statistic $S(t)$ corresponding to $S_n$ is
\begin{equation}
\label{ctstatistic}
S(t) = S_{N(t)} = G_{N(t)}(X_1, \ldots, X_{N(t)}).
\end{equation}
In plain words, the continuous-time statistic corresponds to the statistic of a random number $N(t)$ of random variables $X_i$s.
In order to connect continuous-time statistics and relaxation equations, we consider a special
class of statistics of {\em convolution} type. We will denote these statistics with
the following symbol
\begin{equation}
\label{convstat}
S_n = \underset{i=1}{\overset{n}{\bigoplus}} X_i,
\end{equation}
and we assume the existence of a transform $\mathcal{L}_{\bigoplus}$ such that
\begin{equation}
\label{transstat}
\mathcal{L}_{\bigoplus} (F_{S_n} (u)) (w) = [\mathcal{L}_{\bigoplus} (F_{X_1} (u)) (w)]^n.
\end{equation}
Let us now consider a continuous-time statistic of convolution type (with $N(t)$ independent from the $X_i$s)
\begin{equation}
\label{contconvstat}
S(t) = S_{N(t)} = \underset{i=1}{\overset{N(t)}{\bigoplus}} X_i,
\end{equation}
and let us compute its cumulative
distribution function. We have
\begin{equation}
\label{cumdistrfunct}
F_{S(t)} (u) = \mathbb{P}(S(t) \leq u) = \sum_{n=0}^\infty F_{S_n} (u) \mathbb{P}(N(t) = n).
\end{equation}
Let $Q(w,s)$ denote the Laplace-$\mathcal{L}_{\bigoplus}$ transform of $F_{S(t)} (u)$
\begin{equation}
\label{doubletrans}
Q(w,s) = \mathcal{L} \mathcal{L}_{\bigoplus} (F_{S(t)} (u)) (w,s).
\end{equation}
We have that (see Appendix)
\begin{equation}
\label{MWgen}
Q(w,s) = \mathcal{L} (\bar{F}_{J_1} (t)) (s) \frac{1}{1 - \mathcal{L} (f_{J_1} (t)) (s) \mathcal{L}_{\bigoplus} (F_{X_1} (u)) (w) }.
\end{equation}
Now, following \cite{mainardi00} (see Appendix for details), we can invert the Laplace transform in \eqref{MWgen} to get
\begin{equation}
\mathcal{Q}(w,t) = \mathcal{L}_{\bigoplus} (F_{S(t)} (u)) (w) = \mathcal{L}^{-1} (Q(w,s)) (t),
\end{equation}
$\mathcal{Q}(w,t)$ is the solution of the
Cauchy problem (with initial condition $\mathcal{Q}(w, t=0) = 1$) for the following pseudo-differential relaxation equation
\begin{equation}
\label{relax}
\int_0^t \Phi (t-t') \frac{\partial \mathcal{Q}(w,t')}{\partial t'} \, dt' = - (1 - \mathcal{L}_{\bigoplus} (F_{X_1} (u)) (w)) \mathcal{Q}(w,t),
\end{equation}
where $\mathcal{L} (\Phi(t)) (s)$ is given by eq. \eqref{memorykernel}.

\section{Examples}

To show that the theory developed above is not void, it is possible to present some examples. 

\subsection{Sum statistics}

The first example is the {\em sum statistic} for independent and identically distributed random variables.
\begin{equation}
\label{example1}
S^{(1)}_n = \sum_{i=1}^n X_i;
\end{equation}
the corresponding continuous-time sum statistic is an uncoupled continuous-time random walk:
\begin{equation}
\label{example1bis}
S^{(1)} (t) = \sum_{i=1}^{N(t)} X_i.
\end{equation}
where we take $N(t)$ to be the Poisson process;
in this case, $\bigoplus$ is the usual convolution and the operator $\mathcal{L}_{\bigoplus}$ coincides with the usual Laplace transform $\mathcal{L}$. As we have exponentially distributed $J_i$s, we recall that the kernel $\Phi(t)$ in \eqref{relax} coincides with Dirac's delta $\delta (t)$ and eq. \eqref{relax} becomes an ordinary relaxation equation
\begin{equation}
\label{deltakernel}
\frac{\partial \mathcal{Q}^{(1)}(w,t)}{\partial t} = - (1 - \mathcal{L} (F_{X_1} (u)) (w)) \mathcal{Q}^{(1)}(w,t).
\end{equation}
The solution of the Cauchy problem for the above relaxation equation is
\begin{equation}
\label{soldeltakernel}
\mathcal{Q}^{(1)}(w,t) = \exp( - (1 - \mathcal{L} (F_{X_1} (u)) (w)) t)
\end{equation}
leading, upon inversion of the second transform, to
\begin{equation}
\label{solutionsum}
F_{S^{(1)}(t)} (u) = \exp(-t) \sum_{n=0}^\infty F_{X_1} ^{\star n}(u) \frac{t^n}{n!},
\end{equation}
where $F_{X_1} ^{\star n}(u)$ denotes the $n$-fold convolution and $F_{X_1} ^{\star 0}(u) = \delta(u)$.

\subsection{Maximum statistics}

As a second example, we consider the
{\em maximum statistic}
\begin{equation}
\label{example2}
S^{(2)}_n = \max(X_1, \ldots, X_n);
\end{equation}
this time we assume that the interarrival times $J_i$ follow a Mittag-Leffler distribution of order $0<\alpha <1$; this is characterized by the following complementary cumulative distribution function
\begin{equation}
\label{ML}
\bar{F}_{J_1} (t) = E_{\alpha} (- t^\alpha),
\end{equation}
where the one-parameter Mittag-Leffler function $E_\alpha (z)$ is defined as
\begin{equation}
\label{MLF}
E_\alpha (z) = \sum_{n=0}^\infty \frac{z^n}{\Gamma(\alpha n + 1)}.
\end{equation}
The corresponding continuous-time maximum statistics is
\begin{equation}
\label{example2bis}
S^{(2)} (t) = \max(X_1, \ldots, X_{N_\alpha (t)}),
\end{equation}
where $N_\alpha (t)$ is the fractional Poisson process of renewal type \cite{mainardi04} and, as before, we assume the independence of $N(t)$ from the $X_i$s; here, $\bigoplus$ is the usual product and the
operator $\mathcal{L}_{\bigoplus}$ is the identity. The kernel is
\begin{equation}
\label{caputo}
\Phi (t) = \frac{t^{-\alpha}}{\Gamma(1-\alpha)},
\end{equation}
and the non-local relaxation equation \eqref{relax} becomes
\begin{equation}
\label{solMLkernel}
\frac{\partial^\alpha \mathcal{Q}^{(2)}(w,t) }{\partial t^\alpha} = -(1 - F_{X_1} (w)) \mathcal{Q}^{(2)}(w,t),
\end{equation}
where $\partial^\alpha/\partial t^\alpha$ is the Caputo derivative (see Appendix). The solution of \eqref{solMLkernel} is
\begin{equation}
\label{solutionmax}
F_{S^{(2)}(t)} (w) = \mathcal{Q}^{(2)}(w,t) = E_\alpha ( -(1-F_{X_1} (w)) t^\alpha).
\end{equation}
If $\alpha =1$ and $X_1$ is exponentially distributed, eq. \eqref{solutionmax} reduces to the well-known Gumbel distribution \cite{gumbel35}
\begin{equation}
\label{gumbel}
F_{S^{(2)}(t)} (w) = \exp ( - \exp(- w) t).
\end{equation}
Incidentally, this result has potential applications in geophysics where power-law distributed interarrival times are often observed between extreme events. This was presented earlier in a master thesis \cite{pissarello10}. A general discussion of this problem can be found in \cite{benson07}.



\section{Summary and conclusions}

In this paper, the relation between generalized anomalous relaxation equations and semi-Markov processes is explored in some specific cases. Explicit evolution equations are given for transforms of the cumulative distribution function of continuous-time statistics of convolution-type. Two specific examples are worked out in detail: the sum statistic and the maximum statistic. The case of the sum statistics coincides with the continuous-time random walk. For the maximum statistic, in the presence of power-law interarrival times following the Mittag-Leffler distribution, the theory leads to an explicit analytic form for the cumulative distribution that was not published before. It is a fractional generalization of the well-known Gumbel distribution and it is given in eq. \eqref{solutionmax}.

The theory outlined in Section III and yielding equations such as \eqref{CTRWevolution}, \eqref{Relax3} and \eqref{relax} is leading to interesting developments presented e.g. in \cite{meerschaert15,georgiou15}. Essentially, the idea is that a random time change $N(t)$, where $N(t)$ is a counting renewal process, in a Markov chain leads to a generalized relaxation equation for relevant probabilities (or characteristic functions) whose solution is given in terms of the complementary cumulative distribution of the inter-event duration. We are currently working on mixing properties and stability of these processes. Moreover, this is quite a rich class of processes and there is virtually no limit to modelling, extensions and generalizations.

\section*{Acknowledgements}
Inspiring discussion with Mark M. Meerschaert and Hans-Peter Scheffler is gratefully acknowledged. This took place during the {\em Workshop on Dependence, Stability, and Extremes} held at the Fields Institute, Toronto, Canada, on May 2-6, 2016 \cite{fields16}. The idea of writing this paper originated at that meeting and was later corroborated by participation to the workshop on {\em Ergodicity breaking and anomalous dynamics} at the University of Warwick, UK on 10-12 August 2016 \cite{warwick16}. A preliminary version of this paper was presented at the {\em Workshop on Future Directions in Fractional Calculus Research and Applications} held at Michigan State University on 17--21 October 2016 \cite{scalas16}.

\section*{Appendix}

Start from eq. \eqref{cumdistrfunct} that can be derived from purely probabilistic considerations:
$$
F_{S(t)} (u) = \mathbb{P}(S(t) \leq u) = \sum_{n=0}^\infty F_{S_n} (u) \mathbb{P}(N(t) = n).
\eqno(A.1)
$$
Use the fact that 
$$
\mathcal{L}(\mathbb{P}(N(t)=n))(s) = \mathcal{L}(\bar{F}_{J_1} (t))(s) \, [ \mathcal{L}(f_{J_1} (t))(s)]^n
\eqno(A.2)
$$
and use this and \eqref{transstat} in \eqref{doubletrans} to get
$$ Q(w,s) =
\mathcal{L}(\bar{F}_{J_1} (t))(s) \sum_{n=0}^\infty [\mathcal{L}(f_{J_1} (t))(s) \mathcal{L}_{\bigoplus} (F_{X_1} (u)) (w)]^n.
\eqno(A.3)
$$
Now use the sum of the geometric series to get to \eqref{MWgen}, a generalization of Montroll-Weiss equation.
Remember that $f_{J_1} (t) = - d\bar{F}_{J_1} (t)/dt$ and therefore
$$
\mathcal{L}(\bar{F}_{J_1} (t))(s) = \frac{1 - \mathcal{L}(f_{J_1} (t))(s)}{s}.
\eqno(A.4)
$$
Eq. \eqref{MWgen} can be rearranged as follows \cite{mainardi00}
$$
\frac{1 - \mathcal{L}(f_{J_1} (t))(s)}{s \mathcal{L}(f_{J_1} (t))(s)} (s Q(w,s) - 1) = 
-(1-\mathcal{L}_{\bigoplus} (F_{X_1} (u)) (w))Q(w,s) \nonumber \;\;\;
\eqno(A.5)
$$
and inverting this with respect to the Laplace transform yields eq. \eqref{relax} with
$$
\mathcal{L}(\Phi(t))(s) = \frac{1 - \mathcal{L}(f_{J_1} (t))(s)}{s \mathcal{L}(f_{J_1} (t))(s)} =
\frac{\mathcal{L}(\bar{F}_{J_1} (t))(s) }{\mathcal{L}(f_{J_1} (t))(s)}
\eqno(A.6)
$$
and initial condition $\mathcal{Q}(w,t=0)=1$.

For the Caputo derivative in the second example, replace the kernel given in eq. \eqref{caputo} in the non-local term of eq. \eqref{relax} to get
$$
\frac{1}{\Gamma(1-\alpha)} \int_0^t (t-t')^{-\alpha} \frac{\partial \mathcal{Q}^{(2)}(w,t')}{\partial t'} \, dt'.
$$
This is indeed the Caputo derivative of $\mathcal{Q}^{(2)}(w,t)$ \cite{mainardi00,baleanu16}.

\end{document}